\documentclass[]{elsarticle}
\usepackage{amsmath}%
\usepackage{amsfonts}%
\usepackage{amssymb}%
\usepackage{graphicx}
\newtheorem{theorem}{Theorem}

\numberwithin{equation}{section}

\newproof{proof}{Proof}

\begin{document}

\title{Strong Convergence in Posets}
\author{Amir Ban}
\ead{amirban@netvision.net.il}
\address{Center for the Study of Rationality, Hebrew University, Jerusalem, Israel}
\author{Nati Linial}
\ead{nati@cs.huji.ac.il} 
\address{School of Computer Science and Engineering, Hebrew University, Jerusalem, Israel}
\date{}

\begin{abstract}

We consider the following solitaire game whose rules are reminiscent of the children's game of leapfrog. The game is played on a poset $(P,\prec)$ with $n$ elements. The player is handed an arbitrary permutation
$\pi=(x_1,x_2,\ldots,x_n)$ of the elements in $P$. At each round
an element may ``skip over'' a smaller element preceding it, i.e. swap positions with it. For example, if $x_i \prec x_{i+1}$,
then it is allowed to move from $\pi$ to the permutation
$(x_1,x_2,\ldots,x_{i-1},x_{i+1},x_i,x_{i+2},\ldots,x_n)$ of $P$'s elements. The player is to carry out such steps
as long as such swaps are possible. When there are several consecutive
pairs of elements that satisfy this condition, the player can choose which pair to swap next. Does the player's choice of
swaps matter for the final permutation or is it uniquely determined by the initial order of $P$'s elements?
The reader may guess correctly that the latter proposition is correct. What may be more
surprising, perhaps, is that this question is not trivial. The proof works by constructing an appropriate
system of invariants.

\end{abstract}

\maketitle

\section{Introduction}

Let $(P,\prec)$ be a poset with $n$ elements. 
We say that the elements $x, y \in P$
are \textbf{comparable} if either $x \prec y$ or $y \prec x$.
Otherwise we say they are incomparable and write
$x \parallel y$. Let $\pi=(x_1,x_2,\ldots,x_n)$ be some permutation of $P$'s elements. A
\textbf{swap} changes this permutation to
$(x_1,x_2,\ldots,x_{i-1},x_{i+1},x_i,x_{i+2},\ldots,x_n)$ for some
index $i$. This swap is \textbf{permissible} if $x_i \prec
x_{i+1}$. A permutation of $P$'s elements is called {\em terminal} if no swap is permissible.
We say that $\pi$ {\em converges} to a terminal permutation $\sigma$ of $P$'s elements
if it is possible to move from $\pi$ to $\sigma$ through some sequence of permissible swaps. 
Clearly, every sequence of permissible swaps is finite, since no two
elements can be swapped more than once, and therefore every permutation converges. The main result in this note states that the convergence is {\em strong}, in the sense that the permutation converges to the same terminal permutation regardless of the order of swaps:

\begin{theorem}
Every permutation $\pi$ of the elements of a finite poset
$(P,\prec)$ converges to exactly one terminal permutation.
\label{poconverge}
\end{theorem}

The question arose in the analysis of certain dynamic economic systems~\cite{Ban} in which agents compete against each other. An agent's standing in the process is quantified by a positive real number that captures his current {\em reputation}. As the process unfolds, reputations go up and down, but the absolute scores do not matter and the only issue is the agents' ranking by reputation. It can be shown that in the setting of~\cite{Ban} the possible order changes are governed by an (unknown, but fixed) partial order on the agents. However, the actual timing and dynamics at which agents move up and down in the reputation ranking is controlled by a complex set of stochastic rules. The question came up, then, what can be said about the players' ranking when the system reaches a steady state. As our result shows, the steady state of the dynamical system is determined solely by its initial state. We hope that the result may be useful in the study of dynamical systems in the physical and biological fields as well.

Elementary though this question is, it seems to be new. We note, however, that strong convergence has been studied before in combinatorics. Several beautiful theorems establish strong convergence phenomena in various ``reflection games''. The common theme is this: A graph $G=(V,E)$ is given along with a function $w: V \rightarrow \mathbb Z$. In addition, there is a rule that allows certain local modifications to $w$. Such local changes can take place only as long as $w$ takes some negative values and the process terminates when $w \ge 0$. Several such games were analyzed in the literature~\cite{Alon, Mozes, Eriksson, Winkler, Stanley}. All these papers show that even though there generally exist a number of possible allowed steps, every sequence of allowed steps terminates and the terminal position is independent of the steps chosen. In some cases it is even shown that the number of steps till termination depends only on the initial position. Despite the apparent similarity, we do not see a reduction between our results and these theorems.

\section{The Proof}

\begin{proof}
Let $\pi$ converge to some terminal permutation $\tau$.
The uniqueness of $\tau$ is proved by providing a criterion, depending only on $\pi$, as to
which pairs of elements appear in the same order in $\pi$ and $\tau$ and which are reversed.

An $(x,y)$-{\em link} in $\pi$ is a sequence $x=z_1,
z_2, \ldots, z_k = y$ that appear in this order (not necessarily
consecutively) in $\pi$ such that $z_{\alpha} \parallel z_{\alpha+1}$
for every $\alpha \in [k-1]$. 

Clearly, if $x \parallel y$ no permissible swap can change the relative order of $x$ and $y$.
Consequently:

\begin{itemize}
\item No sequence of permissible swaps can change the relative order of $x$ and $y$ if an $(x,y)$-{\em link} exists.
\item No sequence of permissible swaps can create or eliminate an $(x,y)$-{\em link}.
\end{itemize}

We say that $(x,y)$ is a {\em critical pair} in $\pi$ if (i) $x \prec y$, (ii) $x$
precedes $y$ in $\pi$ and, (iii) there is no $(x,y)$-{\em link} in $\pi$.

We now assert and prove the criterion for whether any two elements $x$ and $y$ in $\pi$, with $x$ preceding $y$,
preserve or reverse their relative order in a terminal permutation:

\begin{enumerate}
\item If $y \prec x$, the order is preserved.
\item If there exists an $(x,y)$-{\em link}, the order is preserved.
\item Otherwise, i.e. if $(x,y)$ is a critical pair, the order is reversed.
\end{enumerate}

The first element of the criterion is trivial and the second has already been dealt with.
It remains to show the third and last element: Since an $(x,y)$-{\em link} cannot be created
or eliminated by permissible swaps, an equivalent statement to this claim is that a permutation $\tau$ with a
critical pair cannot be terminal. We prove this by induction on the number of elements in $\tau$
separating $x$ and $y$:

At the base of induction, if $x$ and $y$ are neighbor elements, the assertion is true since as $x \prec y$ the permutation
is not terminal. Now let $k$ be the number of elements separating $x$ and $y$, and the induction hypothesis is 
that if the number of elements separating a pair is less than $k$ it cannot be critical.

Let $z$ be an element between $x$ in $y$ in $\tau$. Assume $x \prec z$. Then by the induction hypothesis there
exists an $(x,z)$-{\em link}. Now consider the relation between $z$ and $y$: $z \parallel y$ is impossible,
because then $y$ could be concatenated to the $(x,z)$-{\em link} to form an $(x,y)$-{\em link}, contrary to
the assumption that $(x,y)$ is a critical pair. Similarly, $z \prec y$ would by the induction hypothesis
prove the existence of a $(z,y)$-{\em link}, but this is impossible as it could be concatenated to the
$(x,z)$-{\em link} to form an $(x,y)$-{\em link}. This leaves $y \prec z$ as the only possibility.

In summary $x \prec z \Rightarrow y \prec z$.

By similar reasoning $z \prec y \Rightarrow z \prec x$.

Furthermore, the possibility $x \parallel z$ together with $z \parallel y$ can be
dismissed as constituting an $(x,y)$-{\em link}, leaving just two possible scenarios
satisfied by each $z$ between $x$ and $y$:

\begin{itemize}
\item $z$ is ``small'', i.e. $z \prec x$ and $z \prec y$.
\item $z$ is ``large'', i.e. $x \prec z$ and $y \prec z$.
\end{itemize}

Since $\tau$ is terminal, $x$'s immediate neighbor must be ``small'', and $y$'s immediate neighbor
must be ``large''. Between these two, there must exist two consecutive elements $z_1, z_2$ such that
$z_1$ is ``small'' and $z_2$ is ``large''. But this leads to a contradiction as $z_1 \prec x \prec z_2 \Rightarrow z_1 \prec z_2$, which
implies that $\tau$  is not terminal. Therefore $(x,y)$ cannot be critical, completing the proof by induction that a terminal
permutation cannot have a critical pair.

This completes the proof for the uniqueness of the terminal permutation.
\end{proof}

\section{Remarks}

As the proof shows, the relative final order of every pair of elements $x$ and $y$ is determined by the existence
of a link that connects them. If a link exists, then $x, y$ maintain their initial relative order in the
final permutation. If no such link exists, their final relative order agrees with their mutual order relation.

The proof also demonstrates that the number of swaps to reach termination
is uniquely defined: Between the initial and terminal permutation, sum the displacements of elements that moved forward (alternatively sum the displacements of backward moving elements). 

The result shows as well how to efficiently determine the terminal permutation given how the elements are ordered initially:
Perform an arbitrary swap until no more swaps are possible.

One of the referees asked us whether the result extends from posets to directed acyclic graphs. At least at its simplest form the analogous claim is incorrect. Let $D$ be the directed path on 3 vertices $A \rightarrow B \rightarrow C$. The permutation $ABC$ converges to $BAC$ if $A,B$ are swapped, but to $ACB$ if $B,C$ are swapped.

It is of interest to understand how many terminal permutations various posets have. In particular how is this parameter distributed over all $n$-element posets?

\end{document}